\documentclass{article}
\usepackage{amssymb}
\usepackage{amsmath}

\input{tcilatex}
\begin{document}

\label{inent copy(2)}\bigskip \label{inent copy(3)}

\begin{center}
{\Huge Control of Dams When the Input Is a L\'{e}vy Type Process}

\bigskip

Mohamed Abdel-Hameed

Department of Statistics

College of Business and Economics

UAE University
\end{center}

\bigskip

\bigskip {\LARGE Abstract}

Zuckermann [10] considers the problem of optimal control of a finite dam
using $P_{\lambda ,0}^{M}$ policies, assuming that the input process is
Wiener with drift term \ \ \ \ $\mu \geq 0$. Lam and Lou [7\ ] treat the
case where the input is a Wiener process with a reflecting boundary at its
infimum, with drift term $\mu \geq 0$, using the long-run average and total
discounted cost criteria. Attia [3] obtains results similar to those of Lam
and Lou, through simpler and more direct methods. Bae \textit{et al}. [5]
treat the long-run average cost case when the input process is a compound
Poisson process with a negative drift. In this paper we unify and extend the
results of these authors.



\medskip \noindent {\LARGE 1}. {\LARGE Introduction and summary} \bigskip 
\newline
\indent\ Lam and Lou [7]\ consider the control of a finite dam, with
capacity $V>0$,$\ $ where the water input is a Wiener process, using $%
P_{\lambda ,\tau }^{M}\;$policies. In these policies the water release rate
is assumed to be zero until the water reaches level $\lambda >0.\ $As soon
as this happens the water is released at rate $M>0\ $until the water content
reaches level $\tau >0,$ $\lambda >\tau $.$\ $They use the total discounted
as well as the long-run average cost criteria. Attia [3] obtains the same
results of Lam and Lou, using simpler methods. Lee and Ahn [8] consider the
long-run average cost case, for the $P_{\lambda ,0}^{M}\ $policy, when the
water input is a compound Poisson process. Abdel-Hameed [1] treats the case
where the water input is a compound Poisson process with a positive drift.
He obtains the total discounted as well as the long-run average costs. Bae 
\textit{et al. }[4\ ] consider the $P_{\lambda ,0}^{M}\ $policy in assessing
the workload of an M/G/1 queuing system. Bae \textit{et al. }[5] consider
the log-run average cost for $P_{\lambda ,\tau }^{M}\ $policy in a finite
dam, when the input process is \ a compound Poisson process, with a negative
drift term. At any time, the release rate can be increased from $0$ to $M$
with a starting cost$\;K_{1}M$, or decreased from$\;M$\ to zero with a
closing cost $K_{2}M$.\ Moreover, for each unit of output, a reward $R\;$is
received.\ Furthermore, there is a penalty cost which is a bounded
measurable function on the state space of the content process. In this paper
we treat the more general cases where the input process is assumed to be a
spectrally positive L\'{e}vy or a spectrally positive L\'{e}vy process
reflected at its infimum.

For any process $Y=\{Y_{t},t\geq 0\}\ $with state space $E$,$\;$any Borel
set \ $A$ $\subset E\ $and any functional $f$, $E_{y}(f)\;$denotes the
expectation of\ $f$\ conditional on $Y_{0}=y$,$\;P_{y}(A)\;$denotes the
corresponding probability measure and $\mathbf{I}_{A}(\;)\;$is the indicator
function of the set $A$. In the sequel we will write indifferently $P_{0}\ $%
or $P\ \ $and $E\ _{0}$\ or $E$. Throughout, we let $R=(-\infty ,\infty )$, $%
R_{+}=[0,\infty ),$\ $N=\{1,2,...\}$ and $N_{+}=\{0,1,...\}$. For $x,y\in R$%
, we define $x\vee y=x\max y$ $\ $and $x\wedge y=x\min y$. For every $t\geq
0 $, we define $\underline{Y_{t}}=\underset{0\leq s\leq t}{\inf }(Y_{s}$,$%
\wedge 0),\ \overset{-}{Y}_{t}=\underset{0\leq s\leq t}{\sup }(Y_{s},\vee 0)$%
.

We will use the term "increasing" to mean "non-decreasing" throughout this
paper.

In Section\ 2$,$\ we discuss the cost functionals. In Section 3 we define
the input processes and discuss their properties. In Section 4 we obtain
formulas needed for computing the cost functionals using the total
discounted as well as the long-run average cost cases. In section 5 we
discuss the special cases where the input process is a Gaussian process, a
Gaussian process reflected at its infimum and a spectrally positive L\'{e}vy
process of bounded variation.

\bigskip {\LARGE 2.} {\Large The cost functionals}{\LARGE \ \ \ \ \ \ \ \ \
\ \ \ }

{\LARGE \bigskip }\noindent For each $t\in R_{+}$, let $Z_{t}\ $ be the dam
content at time $t\ $,$\ Z=\{Z_{t},t\in R_{+}\}$.\ We define the following
sequence of stopping times : 
\begin{eqnarray*}
\overset{\wedge }{T}_{0} &=&\inf \{t\geq 0:Z_{t}\geq \lambda \},\;\;\;\;\;%
\overset{\ast }{T}_{0}=\inf \{t\geq \overset{\wedge }{T}_{0}:Z_{t}=\tau \},
\\[0.01in]
\overset{\wedge }{T}_{n} &=&\inf \{t\geq \overset{\ast }{T}_{n-1}:Z_{t}\geq
\lambda \},\;\overset{\ast }{T}_{n}=\inf \{t\geq \overset{\wedge }{T}%
_{n}:Z_{t}=\tau \}\text{, }n=1,2,...\;(2.1)
\end{eqnarray*}%
It follows that the process $Z\ \ $is\ a delayed regenerative process with
regeneration points $\{\overset{\ast }{T}_{n},n=0,1,...\}$.$\ $The
regeneration cycle is defined to be the time between successive regeneration
points. \ During a given cycle, the release rate is either $0$ or $M$. \
When the release rate is zero, the process $Z$ \ is either a spectrally
positive L\'{e}vy process (denoted by $I$) or a spectrally positive L\'{e}vy
process reflected at its infimum (denoted by $Y$), and remains so till the
water reaches level$\ \lambda $; from then until it reaches level $\tau $
the content process behaves like the process $\overset{(M)}{I}=I-M$\ \
reflected at $V\ $, we denote this process by $\overset{\ast }{I}$.$\ $It
follows that, for each $t\geq 0,$\bigskip 
\begin{equation*}
\overset{\ast }{I}_{t}=\overset{(M)}{I}_{t}-\underset{0\leq s\leq t}{\sup }((%
\overset{(M)}{I}_{t}-V)\vee 0)\text{. \ }\ (2.2)
\end{equation*}

When the release rate is $0$, the dam is maintained at a net maintenance
cost rate $g$, where $g\ $is a bounded measurable function on $(l,\lambda )$%
, where $l\ \ $is the lower bound of the state space of the process $I$.
Furthermore, maintenance of the dam when the release rate is $M$ is done at
a cost rate $g^{\ast }$, where $g^{\ast }\ $is a bounded measurable function
on $(\tau ,V].$

For $x\,<\lambda $ and $\alpha \in R_{+}$, the discounted cost during the
interval $[0,\overset{\wedge }{T}_{0})$, denoted by $C_{\alpha }(x,\lambda
,0)$, is given as follows

\begin{equation*}
C_{\alpha }(x,\lambda ,0)=E_{x}\int_{0}^{\overset{\wedge }{T}_{0}}e^{-\alpha
t}g(Z_{t})dt\text{.\ \ \ }(2.3)
\end{equation*}

\bigskip Define, $\overset{\ast }{T^{-}}_{\tau }=\overset{\ast }{T}_{0}-%
\overset{\wedge }{T}_{0}$, when the release rate is $M$, starting at $x\geq
\tau $, the expected discounted cost in the interval $[0,\overset{\ast }{%
T^{-}}_{\tau })$, denoted by $C_{\alpha }(x,\tau ,M)$ is given as follows

\begin{equation*}
C_{\alpha }(x,\tau ,M)=E_{x}\int_{0}^{\overset{\ast }{T^{-}}_{\tau
}}e^{-\alpha t}g^{\ast }(Z_{t})dt\text{. \ \ }(2.4)
\end{equation*}

We now discuss the computations of the cost functionals using the total
discounted cost as well as the long-run average cost criteria. \ Let $%
C_{\alpha }(x)\;$be the expected cost during the first cycle, $[0,\overset{%
\ast }{T}_{0})$,$\;$when $Z_{0}=x.\;$ From the definition of the $P_{\lambda
,\tau }^{M}\;$policy, it follows that for $\lambda<x<V$

\begin{equation*}
C_{\alpha }(x)=M\{K_{1}-RE_{_{x}}\int_{0}^{\overset{\ast }{T^{-}}_{\tau
}}e^{-\alpha t}dt\}+\;C_{g^{\ast }}^{\alpha }(M,x,\tau )\text{, }\lambda <x%
\text{\ }<V\ \ \ (2.5)
\end{equation*}%
and for $x\leq \lambda $

\begin{eqnarray*}
C_{\alpha }(x) &=&M\{K_{2}+\ K_{1}\ E_{x}[e^{-\alpha \overset{\wedge }{T}%
_{0}}]-\frac{R}{\alpha }\{E_{_{x}}[e^{-\alpha \overset{\wedge }{T}_{0}}]\
-E_{_{x}}[e^{-\alpha \overset{\ast }{T}_{0}}]\ \} \\
&&+C_{\alpha }(x,\lambda ,0)+E_{_{x}}[e^{-\alpha \overset{\wedge }{T}%
_{0}}C_{\alpha }((Z_{\overset{\wedge }{T}_{0}}\wedge V),\tau ,M)]\text{.}\ \
\ \ \ \ (2.6)\ \ \ 
\end{eqnarray*}

\ Let $C_{\alpha }(\lambda ,\tau )$ and $C(\lambda ,\tau )\ $denote the
total discounted cost and the long-run average cost, respectively. By
modifying the result in Abdel-Hameed [1], it follows that%
\begin{equation*}
C_{\alpha }(\lambda ,\tau )=C_{\alpha }(x)+\frac{E_{x}[\exp (-\alpha \overset%
{\ast }{T}_{0})]C_{\alpha }(\tau )}{1-E_{\tau }[\exp (-\alpha \overset{\ast }%
{T}_{0})]},\ (2.7)
\end{equation*}%
and

\begin{equation*}
C(\lambda ,\tau )=\frac{M{\Large (}K+RE_{0}[T_{\lambda -\tau }^{+}]{\Large )}%
+C_{\alpha }(x,\lambda ,0)+E_{_{x}}[C_{\alpha }((Z_{\overset{\wedge }{T}%
_{\lambda }}\wedge V),\tau ,M)]}{E_{\tau }[\overset{\ast }{T}_{0}\ ]}-RM%
\text{.\ }\ \ (2.8)
\end{equation*}%
where $K=$\ $K_{1}+\ K_{2}\ $and $C_{\alpha }(\tau )\ $is the total
discounted cost during the interval $[0,\overset{\ast }{T}_{0})$,\ given
that $Z_{0}=\tau $. \bigskip \newline
\bigskip \noindent {\LARGE 3.~The input processes and their characteristics}

\bigskip In this paper we consider the cases where the input process is a
spectrally positive L\'{e}vy process, and\ a spectrally positive L\'{e}vy
process reflected at its infimum. In the remainder of this section we
describe these processes and discuss some of their characteristics.\ The
reader is referred to [6] for a more detailed discussion of the definitions
and results mentioned in this section.

Definition. A L\'{e}vy process \ $L=$ $\{L_{t},t\geq 0\}$ with$\ $ state
space $R\ $is said to be spectrally positive L\'{e}vy process, it has no
negative jumps.\newline
It follows that, for each $\theta \in R_{+},x\in R$,

\begin{equation*}
E[e^{-\theta \ L\ _{t}}]=e^{t\phi (\theta )},
\end{equation*}%
where

\begin{equation*}
\phi (\theta )=-a\theta +\frac{\theta ^{2}\sigma ^{2}}{2}-\int_{0}^{\infty
}(1-e^{-\theta x}-\theta x\mathbf{I}_{\{x<1\}})\upsilon (dx)\text{. \ \ }%
(3.1)
\end{equation*}%
The term $a\in R$, $\sigma ^{2}\in R_{+}\ $are\ the drift and variance of
the Brownian motion, respectively, and $\upsilon $ is a positive measure on $%
(0,\infty )$ satisfying $\int_{0}^{\infty }(x^{2}\wedge 1)\upsilon
(dx)<\infty $.

The function $\phi \ $is known as the L\'{e}vy exponent, and it is strictly
convex and tends to infinity as $\theta \ $tends to infinity. For $\alpha
\in R_{+}$, we define

\begin{equation*}
\eta (\alpha )=\sup \{\theta :\phi (\theta )=\alpha \}\ \ (3.2)\text{,}
\end{equation*}%
the largest root of the equation $\phi (\theta )=\alpha $. It is seen that
this equation has at most two roots, one of which is the zero root. Note
that,$\ E(L_{1\ })=$ $\int_{1}^{\infty }x\upsilon (dx)+\mu $.$\ $%
Furthermore, $\underset{t\rightarrow \infty }{\lim }L_{t\ }\ =\infty \ $if
and only if $E(L_{1\ })>0,\ $and $\underset{t\rightarrow \infty }{\lim }%
L_{t\ }\ =-\infty \ $if and only if $E(L_{1\ })<0$. Also, if$\ E(L_{1\
})=0,\ $then $\underset{t\rightarrow \infty }{\lim }L_{t\ }$\ does not
exist. Furthermore, $\eta (0)>0,$ if and only if $E(L_{1\ })>0.$

An important case is when the process $L\ $is of bounded variations, i.e., $%
\sigma ^{2}=0$ and $\int_{0}^{\infty }(x\wedge 1)\upsilon (dx)<\infty $. Let

\begin{equation*}
\zeta =-a+\int_{0}^{1}x\upsilon (dx)\text{.}
\end{equation*}%
In this case we can write

\begin{equation*}
\phi (\theta )=\zeta \theta -\int_{0}^{\infty }(1-e^{-\theta x})\upsilon (dx)%
\text{, }(3.3)\text{ }
\end{equation*}%
where necessarily $\zeta $ is strictly positive. \newline

Definition. A L\'{e}vy process \ is said to be spectrally negative if it has
no positive jumps. \ \ \ \ \ \ \ \ \ \ \ \ \ \ \ \ \ \ \ \ \ \ \ \ \ \ \ \ \ 

For any spectrally positive\ \ L\'{e}vy input process\ $L$, we let $\overset{%
\wedge }{L}$\ $\ =-\ L$ throughout this paper.\ It is clear that $L$ is \
spectrally positive if and only if the process $\overset{\wedge }{L}\ \ $is
spectrally negative.

We now introduce tools, which will be central in the rest of this paper.

Definition. For any spectrally positive L\'{e}vy\ process with L\'{e}vy
exponent $\phi \ $and for $\alpha \geq 0$, the $\alpha -$\textit{\ scale
function} $W^{\alpha }:R\twoheadrightarrow R_{+}$, $W^{\alpha }(x)=0\ $for
every $x\,<\,0,$ and on $[0,\infty )$ it is defined as the unique right
continuous increasing function such that

\begin{equation*}
\int_{0}^{\infty }e^{-\beta x}W^{(\alpha )}(x)=\frac{1}{\phi (\beta )-\alpha 
}\text{, }\beta >\eta (\alpha )\ \ (3.4)
\end{equation*}%
We will denote $W^{0\ \text{ }}$by $W$ throughout. For $\alpha \geq 0$, we
have (see (8.24) of[ 6])

\begin{equation*}
W^{(\alpha )}(x)=\sum_{k=0}^{\infty }\alpha ^{k}W^{\ast (k+1)}(x)\text{,}\ \
\ \ (3.5)
\end{equation*}%
where $W^{\ast (k)}\ $is the kth convolution of $W\ \ $with itself. \ \ \ \
\ \ \ \ \ \ \ \ \ \ \ \ \ \ \ \ \ \ \ \ \ \ \ \ \ 

It follows that $W^{(\alpha )}(0+)=0\ $if and only if the process $L$\ is of
unbounded variation. Furthermore, $W^{\alpha }$ right and left
differentiable on $(0,\infty ).\ $By $W_{+}^{(\alpha )^{\prime }}(x)$, we
will denote the right derivative of $W^{(\alpha )}\ $in $x$. \ \ \ \ \ \ \ \
\ \ \ \ \ \ \ \ \ \ \ \ \ \ \ \ \ \ \ \ \ 

The adjoint $\alpha -$\ scale function associated\textit{\ with }$W^{(\alpha
)}\ ($denoted$\ $by $Z^{(\alpha )}$) is defined as follows: \ \ \ \ \ \ \ \
\ \ \ \ \ \ \ \ \ \ \ \ \ \ \ \ \ \ \ \ \ 

Definition. For $\alpha \geq 0$, the \textit{adjoint }$\alpha -$\textit{\
scale }$Z^{(\alpha )}:R_{+}\twoheadrightarrow \lbrack 1,\infty )\ $is
defined as

\begin{equation*}
Z^{(\alpha )}(x)=1+\alpha \int_{0}^{x}W^{\alpha }(x)dx.\ \ \ (3.6)
\end{equation*}%
It follows that as $x\rightarrow \infty $,\ for $\alpha >0$, $,W^{(\alpha
)}(x)\sim \frac{e^{\eta (\alpha )\ x}}{\phi ^{^{\prime }}(\eta (\alpha )\ )}$
and $\frac{Z^{\alpha }(x)}{W^{(\alpha )}(x)}\sim \frac{\alpha }{\eta (\alpha
)}$.

\bigskip {\LARGE 4.} {\LARGE Basic results}

\bigskip

To derive $C_{\alpha }(x,\lambda ,0)\ $, $E_{x}[e^{-\alpha \overset{\wedge }{%
T}_{0}}]$,$\ C(\lambda ,\tau )$, $E_{\tau }[\overset{\wedge }{T}_{0}]\ $we
define the input process $Z\ \ $killed at $\overset{\wedge }{T}_{0}$,
denoted by $X$\thinspace\ such that for every $t\geq 0$

\begin{equation*}
X_{t}=\{Z_{t},t<\overset{\wedge }{T}_{0}\}\text{.}\ (4.1)
\end{equation*}%
It is known that this killed process is a strong Markov process.

For the input process $X$, any Borel set $A\ $contained in the state space
of $X$,$\ t\in R_{+}$,$\ $the probability transition function of this
process is given as follows

\begin{equation*}
P_{t}(x,A)=P_{x}(Z_{t}\in A,t<\overset{\wedge }{T}_{0}\}
\end{equation*}%
and for each $\alpha \in R_{+}$ its $\alpha -$potential is defined as follows

\begin{equation*}
U^{\alpha }(x,A)=\int_{0}^{\infty }P_{t}(x,A)e^{-\alpha t}dt=E_{x}\int_{0}^{%
\overset{\wedge }{T}_{0}}e^{-\alpha t}\boldsymbol{I}_{\{Z_{t}\in A\}}dt\text{%
. \ }(4.2)
\end{equation*}%
We note that for $x<\lambda $

\begin{equation*}
C_{\alpha }(x,\lambda ,0)=U^{\alpha }g(x)\text{. \ }\ (4.3)
\end{equation*}

\ The following lemma will be used extensively throughout this paper.

Lemma. Let $S=\{S_{t},t\geq 0\}\ $ be a strong Markov process. Define, $%
\mathcal{G}=\left\{ \sigma (S_{u}\ ,u\leq t)\right\} _{t\geq 0}$,{\LARGE \ }$%
{\huge \kappa }$ to be any stopping time with respect to $\mathcal{G}\ $and $%
\ U^{\alpha ,S}(x,A)\ $as the $\alpha -$potential of the process $S\ \ $%
killed at ${\huge \kappa }$. Denote the state space of this process by $E.$
Then, for $x\in E$ 
\begin{equation*}
E_{x}[e^{-\alpha {\huge \kappa }}]=1-\alpha U_{\mathbf{I}_{E}}^{\alpha ,S}(x)%
\text{.\ \ \ }(4.4)
\end{equation*}

Proof. From the definition of $U^{\alpha ,S}\ $\ and for any bounded \
measurable function $f$\ \ whose domain is $E$, we have

\begin{eqnarray*}
U^{\alpha ,S}f(x) &=&E_{x}[\int_{0}^{{\huge \kappa }}e^{-\alpha
t}f(S_{t})dt]=\int_{E}f(y)U^{\alpha ,S}(x,dy)\cdot \\
&&
\end{eqnarray*}%
Taking $f\ \ $to be identically equal to one, we have

\begin{equation*}
\frac{1-E_{x}[e^{-\alpha {\huge \kappa }}]}{\alpha }=U_{\mathbf{I}%
_{E}}^{\alpha ,S}(x)\text{.}
\end{equation*}%
The required result is immediate from the last equation above. \ \ $%
\blacksquare $ \ \ \ \ \ \ \ \ \ \ \ \ \ \ \ \ \ \ \ \ \ \ \ \ \ \ \ \ \ 

First we consider the case where, during the period $[0,$ $\overset{\wedge }{%
T}_{0})$,$\ $the process $Z\ \ $is a spectrally positive\ L\'{e}vy process,
denoted by $I$. In this case we note that $\overset{\wedge }{T}_{0}=\inf
\{t\geq 0:\overset{}{I}_{t}\ \geq \lambda \}$ and we will denote it by $%
T_{\lambda }^{+}\ $. Throughout the rest of this paper, for any $a\in R$, we
define $T_{a}^{-}=\inf \{t\geq 0:\overset{}{I}_{t}\ \leq a\}$, $\top
_{a}^{+}=\inf \{t\geq 0:\overset{\wedge }{I}_{t}\ \geq a\}$ and $\top
_{a}^{-}=\inf \{t\geq 0:\overset{\wedge }{I}_{t}\ \leq a\}.$\ \ \ \ \ 

\bigskip

Proposition. For $\alpha \geq 0$,\ $a\leq \lambda $ the \ $\alpha $potential 
$\ $($\overset{(1)}{U^{\alpha }}$)$\ $of the process $I\ \ $killed at\ $%
T=T_{\lambda }^{+}\wedge $\ $T_{a}^{-}\ $is absolutely continuous with
respect to the Lebesgue measure on $[a,\lambda ]\ $and a version of its
density is given by\ 

\begin{equation*}
\overset{(1)}{u^{\alpha }}(x,y)=W^{(\alpha )}(\lambda -x)\frac{W^{(\alpha
)}(y-a)}{W^{(\alpha )}(\lambda -a)}-W^{(\alpha )}(y-x),\ \ x,y\in \lbrack
a,\lambda ]\text{.}\ \ \ \ \ \ (4.5)
\end{equation*}

Proof. For $A\subset \lbrack a,\lambda ]$

\begin{eqnarray*}
\overset{(1)}{U}(x,A) &=&E_{x}\int_{0}^{T}e^{-\alpha t}\mathbf{I}%
_{\{I_{t}\in A\}}dt \\
&=&E_{-x}\int_{0}^{\top _{-\lambda }^{-}\wedge T_{-a}^{+}}e_{{}}^{-\alpha t}%
\mathbf{I}_{\{\overset{\wedge }{I}_{t}\in -A\}}dt \\
&=&E_{\lambda -x}\int_{0}^{\top _{0}^{-}\wedge T_{\lambda
-a}^{+}}e_{{}}^{-\alpha t}\mathbf{I}_{\{\overset{\wedge }{I}_{t}\in \lambda
-A\}}dt \\
&=&\int_{(\lambda -A)}{\Large [}W^{(\alpha )}(\lambda -x)\frac{W^{(\alpha
)}(\lambda -a-y)}{W^{(\alpha )}(\lambda -a)}-W^{(\alpha )}(y-x)]dy\text{,}
\end{eqnarray*}%
where the last equation follows from Theorem 8.7 of [7], this establishes
our assertion. $\ \ \ \ \ \ \ \ \ \ \ \blacksquare $\ \ 

\ \ \ \ \ \ \ \ \ \ \ \ \ \ \ \ \ \ \ \ 

Corollary. For $\alpha \geq 0$ the $\alpha $-potential $(U^{\alpha })\ $of
the process $X$\ \ is absolutely continuous with respect to the Lebesgue
measure on $(-\infty ,\lambda ]\ $and a version of its density is given by

\bigskip

\begin{equation*}
u^{\alpha }(x,y)=W^{(\alpha )}(\lambda -x)e^{-\eta (\alpha )(\lambda
-y)}-W^{(\alpha )}(y-x),\ \ x,y\in (-\infty ,\lambda ]\text{. \ }(4.6)
\end{equation*}

Proof: The proof follows from (4.5) by letting $a\rightarrow -\infty $ and
since,\ for $\alpha \geq 0$, $W^{(\alpha )}(x)\sim \frac{e^{\eta (\alpha )\
x}}{\phi ^{^{\prime }}(\eta (\alpha )\ )}\ $as $x\rightarrow \infty $. \ \ \
\ $\ \blacksquare $\ $\ $

\bigskip

With the help of the last corollary above we are now in a position to find $%
E_{x}[e^{-\alpha T_{\lambda }^{+}}]$ and $E_{x}[T_{\lambda }^{+}]$.

\bigskip Theorem. (i) For $\alpha >0$ and $x<\lambda \ $we$\ $have

\begin{equation*}
E_{x}[e^{-\alpha T_{\lambda }^{+}}]=Z^{(\alpha )}(\lambda -x)-\frac{\alpha }{%
\eta (\alpha )}W^{(\alpha )}(\lambda -x)\text{. }(4.7)
\end{equation*}

\ \ \ \ \ \ \ \ \ \ \ \ (ii) For $x<\lambda \ $we$\ $have%
\begin{eqnarray*}
E_{x}[T_{\lambda }^{+}] &=&\frac{W(\lambda -x)}{\eta (0)}-\overset{-}{W}%
(\lambda -x)\text{, \ }\eta (0)>0\text{\ \ }(4.8) \\
&=&\infty \ \ \ \ \ \ \ \ \ \ \ \ \ \ \ \ \ \ \ \ \ \ \ \ \ \ \ \ \ \ \ \ \
\ ,\eta (0)=0\text{, }
\end{eqnarray*}%
where for every $x\geq 0$,%
\begin{equation*}
\overset{-}{W}(x)=\int_{0}^{x}W(y)dy.\ \ (4.9)
\end{equation*}

Proof. We only prove (i), the proof of (ii)\ is easily obtained from (i)\
and hence is omitted. Let $\overset{}{U^{\alpha }}\ $be the $\alpha $-$\ $%
potential of the process $X$, then

\begin{eqnarray*}
E_{x}[e^{-\alpha T_{\lambda }^{+}}] &=&1-\alpha U_{\mathbf{I}_{(-\infty
,\lambda )}}^{\alpha }(x) \\
&=&1-\alpha \int_{-\infty }^{\lambda }{\large \{}W^{\alpha }(\lambda
-x)e^{-(\lambda -y)\eta (\alpha )}-W^{\alpha }(y-x){\Large \}}dy \\
&=&1+\alpha \int_{x}^{\lambda }W^{\alpha }(y-x)dy-\alpha W^{\alpha }(\lambda
-x)\int_{-\infty }^{\lambda }e^{-(\lambda -y)\eta (\alpha )}dy \\
&=&Z^{(\alpha )}(\lambda -x)-\frac{\alpha }{\eta (\alpha )}W^{(\alpha
)}(\lambda -x)\text{,}
\end{eqnarray*}%
where the first equation follows from (4.4), the second equation follows
from (4.6), the third equation follows since $W^{(\alpha )}(x)=0$, $x\,<0$.\
and the last equation follows from the definition of $Z^{(\alpha )}$. \ \ \
\ $\blacksquare \ \ $

For any Borel set $B\subset R_{+}\times R$, we let $M(B)$ be the Poisson
random measure counting the number of jumps of the process $I$ \ in $B$ with
L\'{e}vy measure $\nu $,$\ $where if $B=[0,t)\times A$, $A\subset R$, then $%
E[M(B)]=t\upsilon (A)$. We need the following to compute the last term in
(2.6). \ \ \ \ \ \ \ \ \ \ \ \ \ \ \ \ \ \ \ \ \ \ \ \ \ \ \ \ \ 

Proposition. For $\alpha \geq 0$ let $\overset{(1)}{u^{\alpha }}(x,y)\ $be
as given in (4.5) and $x\leq \lambda \leq z$, then

\begin{equation*}
E_{_{x}}[e^{-\alpha T_{\lambda }^{+}},I_{T_{\lambda }^{+}}\in dz,\
T_{\lambda }^{+}<T_{a}^{-}]=\int_{a}^{\lambda }\upsilon (dz-y)\overset{(1)}{%
u^{\alpha }}(x,y)dy\ \ \ \ \ (4.10)
\end{equation*}

Proof. For $x<\lambda ,\alpha \geq 0$, $C\subset \lbrack \lambda ,\infty )\ $%
and $D\subset (a,\lambda )$ 
\begin{eqnarray*}
E_{x}[e^{-\alpha T_{\lambda }^{+}}\text{, }X_{T_{\lambda }^{+}} &\in
&C,X_{T_{\lambda }^{+}-}\in D,T_{\lambda }^{+}<T_{a}^{-}] \\
&=&E_{x}{\LARGE [}\int_{[0,\infty )\times (0,\infty ),}e_{{}}^{-\alpha t}%
\mathbf{I}_{\{\overset{-}{X}_{t}-<\lambda ,\underset{-}{X_{t-}},>a,X_{t}-\in
D\}}\mathbf{I}_{\{y\in C-X_{t}-\}}M(dt,dy) \\
&=&E_{x}{\LARGE [}\int_{[0,\infty )}e^{-\alpha t}\mathbf{I}_{\{\overset{-}{X}%
_{t-}<\lambda ,\underset{-}{X_{t-}},>a\},}\mathbf{I}_{\{X_{t}\in D\}}\nu
(C-X_{t})dt] \\
&=&E_{x}{\LARGE [}\int_{[0,\infty )),}e^{-\alpha t}\mathbf{I}_{\{t<T\},}%
\mathbf{I}_{\{X_{t}\in B\}}\nu (C-X_{t})dt)] \\
&=&E_{x}{\LARGE [}\int_{[0,\infty )}e^{-\alpha t}\mathbf{I}_{\{t<T\},}\nu
(C-X_{t})\mathbf{I}_{\{X_{t}\in D\}}dt)] \\
&=&E_{x}{\LARGE [}\int_{[0,\infty )\times D}e^{-\alpha t}\mathbf{I}%
_{\{t<T\},}\nu (A-y)\mathbf{I}_{\{X_{t}\in dy\}}dt] \\
&=&\int_{D}\nu (C-y)\overset{(1)}{u^{\alpha }}(x,y)dy,
\end{eqnarray*}%
where the second equation follows from the \textit{compensation formula\ }%
(Theorem 4.4. of [7])\textit{. }Our assertion is proved by taking $%
D=[a,\lambda ]$. \ \ \ \ \ \ \ \ $\blacksquare $ \ \ \ \ \ \ \ \ 

\ \ \ \ \ \ \ \ \ \ \ \ \ \ \ \ \ \ \ \ \ \ \ \ \ \ \ 

The following corollary gives the formula needed to compute the last term of
(2.6), when the input process is a spectrally positive L\'{e}vy process.

Corollary. Let $u^{\alpha }$be\ as defined in (4.6). For $\alpha \geq 0$ and
for $x\leq \lambda \leq z$,

\begin{equation*}
E_{_{x}}[e^{-\alpha T_{\lambda }^{+}},I_{T_{\lambda }^{+}}\in dz\
]=\int_{-\infty }^{\lambda }\upsilon (dz-y)u^{\alpha }(x,y)dy\text{. }\ \ \
(4.11)
\end{equation*}

Proof. The proof follows immediately from (4.6) and (4.10). by letting $%
a\rightarrow -\infty $.\ \ \ \ \ \ \ \ \ \ \ \ \ \ \ \ \ \ \ \ \ \ \ $%
\blacksquare $ \ 

We now turn our attention to the case where, during the period $[0,$ $%
\overset{\wedge }{T}_{0})$,$\ $the process $Z\ \ $is a spectrally positive\
\ L\'{e}vy process reflected at its infimum, denoted by $Y$. We will denote $%
\overset{\wedge }{T}_{0}\ $by ${\huge \tau }_{_{\lambda }}$, in this case.

\bigskip The following proposition gives the $\alpha -$potential of the
process $X\ $ defined in (4.1).\ 

Proposition. Assume that during the period $[0,$ $\overset{\wedge }{T}_{0})$,%
$\ $the process $Z\ \ $is a spectrally positive\ L\'{e}vy process reflected
at its infimum. Denote the\ $\ \alpha $-$\ $potential of the process $X\ $\
by $\overset{(2)}{U^{\alpha }}$. Then for any $x,y\in \lbrack 0,\lambda )$,

\ \ \ \ \ \ \ \ \ \ \ \ \ \ \ \ \ \ \ \ \ \ \ \ \ \ \ \ \ \ \ \ \ \ \ \ \ \
\ \ \ \ \ \ \ \ \ \ \ \ \ \ \ \ \ \ \ \ \ \ \ \ \ \ \ \ \ \ \ \ \ \ \ \ \ \
\ \ \ \ \ \newline
\begin{equation*}
\overset{(2)}{U^{\alpha }}(x,dy)=\frac{W^{(\alpha )}(\lambda -x)W^{(\alpha
)}(dy)}{W_{+}^{(\alpha )^{\prime }}(\lambda )}-W^{(\alpha )}(y-x)dy\text{, \
\ \ \ }(4.12)
\end{equation*}%
where for $x,y\in \lbrack 0,\lambda )$, $W^{(\alpha )}(dy)=W^{(\alpha
)}(0)\delta _{0}(dy)+W_{+}^{(\alpha )^{\prime }}(y)dy$, and $\delta _{0}\ $%
is the delta measure in zero.

\ \ \ \ Proof. Note that for each $t\geq 0,$%
\begin{eqnarray*}
Y_{t} &=&I_{t}-\underline{I_{t}}\ \ \ \ \ \ (4.13)\ \  \\
&=&\overset{-}{\overset{\wedge }{I}}_{t}-\overset{\wedge }{I}_{t}\text{,}
\end{eqnarray*}%
where the process $I=\{I_{t},t\geq 0\}\ $is a spectrally positive L\'{e}vy
process. The result follows from part (ii) \ of Theorem 8.11\ of [7], since
the process $\overset{}{\overset{\wedge }{I}}\ $is a $\ $is a spectrally
negative L\'{e}vy process.~\ \ \ $\blacksquare $

\bigskip

The following proposition gives $E_{x}[e^{-\alpha {\huge \tau }_{_{\lambda
}}}]$ and $E_{x}[{\huge \tau }_{_{\lambda }}]$. $\ $

Proposition. (i) For $\alpha \geq 0$ and $x<\lambda \ $we$\ $have%
\begin{equation*}
E_{x}[e^{-\alpha {\huge \tau }_{_{\lambda }}}]=Z^{(\alpha )}(\lambda
-x)-W^{(\alpha )}(\lambda -x)\frac{\alpha W^{(\alpha )}(\lambda )}{%
W_{+}^{(\alpha )^{\prime }}(\lambda )}.\ \text{ }(4.14)
\end{equation*}

\ \ \ \ \ \ \ \ \ \ \ \ \ \ \ \ \ \ (ii) For $x<\lambda \ $we$\ $have%
\begin{equation*}
E_{x}[{\huge \tau }_{_{\lambda }}]=W(\lambda -x)\frac{W(\lambda )}{%
W_{+}^{^{\prime }}(\lambda )}-\overset{-}{W}(\lambda -x)\text{. \ }(4.15)
\end{equation*}

Proof. The proof of part (i) follows from (4.4) and (4.12), in a manner
similar to the proof of (4.7). The proof of part (ii) follows from part (i)
be direct differentiation. We omit both proofs.$\ \ \ \ \ \ \ \blacksquare $

\bigskip

To find a formula analogous to (4.11), for the spectrally positive\ L\'{e}vy
process reflected at its infimum, we first need few definitions. For$\
z>\lambda $ we let

\bigskip 
\begin{eqnarray*}
l_{\alpha }(dz) &=&W^{(\alpha )}(\lambda -x)\int_{0}^{\lambda }W^{(\alpha
}(dy)\upsilon (dz-y)\ \ \ \ \ \ \ \ \ \ \ \ \ \ \ \ \ \ \ (4.16) \\
&&-W_{+}^{(\alpha )^{\prime }}(\lambda )\int_{0}^{\lambda }dyW^{(\alpha
)}(y-x)\upsilon (dz-y)]. \\
L_{\alpha }(z) &=&\int_{(z,\infty )}^{{}}l_{\alpha }(du)\text{.\ \ \ \ \ \ \
\ \ \ \ \ \ \ \ \ \ \ \ \ \ \ \ \ \ \ \ \ \ \ \ \ \ \ \ \ \ \ \ \ \ \ \ \ }%
(4.17) \\
V_{\alpha }(\lambda ) &=&W_{+}^{(\alpha )^{\prime }}(\lambda )Z^{(\alpha
)}(\lambda -x)-\alpha W^{(\alpha )}(\lambda -x)W^{(\alpha )}(\lambda ).\text{%
\ }(4.18)\ \text{\ }
\end{eqnarray*}

The following proposition gives the required formula.

Proposition. (i) For $\alpha \geq 0$ and for $x\leq \lambda <$ $z$,

\begin{equation*}
E_{_{x}}[e^{-\alpha {\huge \tau }_{_{\lambda }}},Y_{{\huge \tau }_{_{\lambda
}}}\in dz]=\frac{l_{\alpha }(dz)}{W_{+}^{(\alpha )^{\prime }}(\lambda )}\ 
\text{,\ }z>\lambda .\ \ \ \ \ (4.19)
\end{equation*}

\ \ \ \ \ \ \ \ \ \ \ \ \ \ \ \ \ (ii) For $\alpha \geq 0$

\begin{equation*}
E_{_{x}}[e^{-\alpha {\huge \tau }_{_{\lambda }}},Y_{{\huge \tau }_{_{\lambda
}}}=\lambda ]=\frac{V_{\alpha }(\lambda )-L_{\alpha }(\lambda )}{%
W_{+}^{(\alpha )^{\prime }}(\lambda )}.\ \ (4.20)\ 
\end{equation*}

Proof. (i) From (4.13), for $x\geq 0$, $Y_{0}\ =x\ \ $if and only if $%
I_{0}=x $ if and only if $\overset{\wedge }{I}_{0}=-x$ . Furthermore, $Y_{%
{\huge \tau }_{_{\lambda }}}=I_{T_{\lambda }^{+}}$ \ almost surely on $%
\{T_{\lambda }^{+}<T_{0}^{-}\}$.\ Therefore

\begin{eqnarray*}
E_{_{x}}[e^{-\alpha {\huge \tau }_{_{\lambda }}},Y_{{\huge \tau }_{_{\lambda
}}} &\in &dz\}]=E_{_{x}}[e^{-\alpha {\huge \tau }_{_{\lambda }}},\ Y_{{\huge %
\tau }_{_{\lambda }}}\in dz,T_{\lambda }^{+}<T_{0}^{-}]+E_{_{x}}[e^{-\alpha 
{\huge \tau }_{_{\lambda }}},Y_{{\huge \tau }_{_{\lambda }}}\in
dz,T_{\lambda }^{+}>T_{0}^{-}] \\
&=&E_{_{x}}[e^{-\alpha T_{\lambda }^{+}}\ ,I_{T_{\lambda }^{+}}\in
dz,T_{\lambda }^{+}<T_{0}^{-}]+E_{x}[e^{-\alpha T_{0}^{-}}\text{,}T_{\lambda
}^{+}>T_{0}^{-}] \\
&&\times E_{_{0}}[e^{-\alpha {\huge \tau }_{_{\lambda }}},Y_{{\huge \tau }%
_{_{\lambda }}}\in dz] \\
&=&E_{_{x}}[e^{-\alpha T_{\lambda }^{+}}\ ,I_{T_{\lambda }^{+}}\in
dz,T_{\lambda }^{+}<T_{0}^{-}]+E_{-x}{\Large [}e^{-\alpha T_{0}^{-}},\top
_{-\lambda }^{-}>\top _{0}^{+}{\Large ]} \\
&&\times E_{_{0}}{\Large [}e^{-\alpha {\huge \tau }_{_{\lambda }}},Y_{{\huge %
\tau }_{_{\lambda }}}\in dz{\large ]} \\
&=&E_{_{x}}[e^{-\alpha T_{\lambda }^{+}}\ ,I_{T_{\lambda }^{+}}\in
dz,T_{\lambda }^{+}<T_{0}^{-}]+E_{\lambda _{-x}}[e^{-\alpha \top _{\lambda
}^{+}},\top _{0}^{-}>\top _{\lambda }^{+}]\ \ \ \ \ \ \ \ \ \ \ \  \\
&&\times E_{_{0}}[e^{-\alpha {\huge \tau }_{_{\lambda }}},Y_{{\huge \tau }%
_{_{\lambda }}}\in dz]\text{{\large , \ \ \ \ \ \ \ \ \ \ \ \ \ \ \ \ \ }}\
(4.21)\text{{\large \ \ \ \ \ \ \ \ \ \ \ \ \ \ \ \ \ \ \ \ }}
\end{eqnarray*}%
where the second equation follows using the strong Markov property, the
third and fourth equations follow from since $\overset{\wedge }{I}=-\overset{%
}{I}$\ and from the definitions of $T_{\lambda }^{+},T_{0}^{-}\ ,\top
_{0}^{-},\top _{\lambda }^{+}$. \ 

Letting $a\rightarrow 0\ $in (4.5) and (4.10), we find that the first term
in the last equation above is equal to $\dint\limits_{0}^{\lambda }\nu (dz-y)%
{\Large [}W^{(\alpha )}(\lambda -x)\frac{W^{(\alpha )}(y)}{W^{(\alpha
)}(\lambda )}-W^{(\alpha )}(y-x){\large ]}dy$. The second term is equal to $%
\frac{W^{(\alpha )}(\lambda -x)}{W^{(\alpha )}(\lambda )}\ $(see (8.8) of
[6]) and the third term is equal to $\frac{h_{\alpha }(dz)}{W_{+}^{(\alpha
)^{\prime }}(\lambda )}$ (this follows from Theorem 4.1 of [9] by letting
the $\beta ,\gamma \rightarrow 0$ and noting that if $\Pi (dz)\ $is the L%
\'{e}vy$\ $measure of the process $\overset{\wedge }{I}$, then (for all $%
z\geq 0$) $\Pi (-\infty ,-z]=\nu \lbrack z,\infty )$.

Our assertion is satisfied by replacing each of the three terms in (4.21)\
by the corresponding value indicated in the last paragraph and after some
algebraic manipulations, which we omit.

\ \ \ \ \ \ \ \ \ (ii) The proof is immediate from (4.14) and (4.19). \ \ \
\ \ \ \ \ \ \ \ \ \ \ \ \ \ $\ \blacksquare $ \ \ 

\bigskip Now we turn our attention to computing $C_{\alpha }(x,M,\tau )$, $%
E_{x}[\exp (\overset{\ast }{T^{-}}_{\tau })]$, and $E_{x}[\overset{\ast }{%
T^{-}}_{\tau }]$. For each $t\geq 0$
\begin{equation*}
\overset{\ast }{X}_{t}=\{\overset{\ast }{I}_{t},t<\overset{\ast }{T^{-}}%
_{\tau }\}\text{.}\ \ \ \ (4.22)
\end{equation*}%
We note that the sample paths of a spectrally positive L\'{e}vy process and
a spectrally positive L\'{e}vy process reflected at its infimum behave the
same way starting at any $x>\tau \ $until they reach level $\tau $, thus $%
\overset{\ast }{X}\ $behaves the same way in both cases. Let $\overset{\ast }%
{U}^{\alpha }$be the $\alpha $potential of the process $\overset{\ast }{X}$.
For each $x\in (\tau ,V]$ 
\begin{equation*}
C_{\alpha }(x,M,\tau )=\overset{\ast }{U}^{\alpha }g^{\ast }(x)\ \ \ (4.23)
\end{equation*}

Let $\overset{(M)}{I}=I-M$\ , as defined in Section 2.1. We note that this
process is a spectrally positive L\'{e}vy process \ with the L\'{e}vy
exponent $\phi _{M}(\theta )=\phi (\theta )+\theta M$,$\ \theta \bigskip
\geq 0$. $\ $We denote its $\alpha -$scale and adjoint $\alpha -$scale
functions by $W_{M}^{(\alpha )}$ and $Z_{M}^{(\alpha )}$, respectively. Note
that $W_{M}^{(\alpha )}\ $is obtained from $W_{{}}^{(\alpha )}$, upon
replacing the term $a$ in (3.1) by $a-M.$

Theorem. For $\alpha >0,$ $\overset{\ast }{U}^{\alpha }\ $is absolutely
continuous with respect to the Lebsegue measure on $(\tau ,V]$, and a
version of its density is given by

\begin{equation*}
\overset{\ast }{u}^{\alpha }(x,y)=\frac{Z_{M}^{(\alpha )}(V-x)W_{M}^{(\alpha
)}(y-\tau )}{Z_{M}^{(\alpha )}(V-\tau )}-W_{M}^{(\alpha )}(y-x)\ .\ \ \
x,y\in (\tau ,V]\text{ \ }(4.24)
\end{equation*}

Proof. For each $t\geq 0$, we define $M\ _{t}=\overset{(M)}{I}_{t}-V$.\ $\ \ 
$For\ any $b\in R$, we define $\gamma _{b}^{+}=\inf \{t\geq 0:\overset{-}{M\ 
}_{t}-M\ _{t}>b\}=\inf \{t\geq 0:\overset{\wedge }{M}_{t}-\underline{\overset%
{\symbol{94}}{M}_{t}}>b\}\ $and $\gamma _{b}^{-}=\inf \{t\geq 0:M\ _{t}-%
\overset{-}{M\ }_{t}<b\}$.\ For and Borel set $A\subseteq (\tau ,V]$ and $%
x\in (\tau ,V]$ we have

\begin{eqnarray*}
P_{x}\{\overset{\ast }{X}_{t} &\in &A\}=P_{x}\{\overset{\ast }{I}_{t}\in
A,t<T_{\tau }^{-}\} \\
&=&P_{x}\{\overset{(M)}{I}_{t}-\underset{s\leq t}{\sup }((\overset{(M)}{I}%
_{s}-V)\vee 0)\in A,t<T_{\tau }^{-}\} \\
&=&P_{x-V}\{M_{t}-\overset{-}{M\ }_{t}\in A-V,t<\tau _{\tau -V}^{-}\} \\
&=&P_{V-x}\{\overset{-}{M\ }_{t}-M_{t}\ \in V-A,t<\tau _{\tau -V}^{-}\} \\
&=&P_{V-x}\{\overset{\wedge }{M}_{t}-\underline{\overset{\symbol{94}}{M}_{t}}%
\ \in V-A,t<\tau _{V-\tau }^{+}\}
\end{eqnarray*}%
Using Theorem 8.111 (ii) of [6], the result follows. \bigskip\ \ \ $%
\blacksquare $

The following theorem gives Laplace transform of the distribution of the
stopping time $\overset{\ast }{T^{-}}_{\tau }\ $and $E_{x}[\overset{\ast }{%
T^{-}}_{\tau }],x\in (\tau ,V]$.

Theorem. (i) Let $x\in \lbrack \tau ,V)$and $\alpha \in R_{+}$. \ Then we
have

\begin{equation*}
E_{x}[e^{-\alpha \overset{\ast }{T^{-}}_{\tau }}]=\frac{Z_{M}^{(\alpha
)}(V-x)}{Z_{M}^{(\alpha )}(V-\tau )}\text{. \ }(4.25)
\end{equation*}

\ \ \ \ \ \ \ \ \ \ \ \ \ \ \ \ \ \ \ \ (ii)\ \ \ For $x\in \lbrack \tau ,V)$%
\begin{equation*}
E_{x}[\overset{\ast }{T^{-}}_{\tau }]=\overset{-}{W_{M}}(V-\tau )-\overset{-}%
{W_{M}}(V-x)\text{,\ \ \ }(4.26)
\end{equation*}%
where, $\overset{-}{W_{M}}(x)=\ \dint\limits_{0}^{x}\overset{}{W_{M}}(y)dy$.

Proof. \ We only prove (i), the proof of (ii) follows easily from (i) and is
omitted. We have 
\begin{eqnarray*}
E_{x}[e^{-\alpha \overset{\ast }{T^{-}}_{\tau }}] &=&1-\alpha \overset{\ast }%
{U}^{\alpha }\mathbf{I}_{(\tau ,V]}(x) \\
&=&1-\alpha \int_{\tau }^{V}\overset{\ast }{u}^{\alpha }(x,dy) \\
&=&1-\alpha \int_{\tau }^{V}{\Large [}\frac{Z_{M}^{(\alpha )}(V-x)W^{\alpha
}(y-\tau )}{Z_{M}^{(\alpha )}(V-\tau )}-W_{M}^{(\alpha )}(y-x){\Large ]}dy\ 
\\
&=&1-\alpha {\Large [}\frac{Z_{M}^{(\alpha )}(V-x)}{Z_{M}^{(\alpha )}(V-\tau
)}{\Large \{}\frac{Z^{\alpha }(V-\tau )-1}{\alpha }{\Large \}}-{\large \{}%
\frac{Z_{M}^{(\alpha )}(V-x)-1}{\alpha }{\large \}}{\LARGE ]} \\
&=&\frac{Z_{M}^{(\alpha )}(V-x)}{Z_{M}^{(\alpha )}(V-\tau )}-Z^{\alpha
}(V-x)+Z^{\alpha }(V-x) \\
&=&\frac{Z_{M}^{(\alpha )}(V-x)}{Z_{M}^{(\alpha )}(V-\tau )}\text{. }
\end{eqnarray*}%
where the third equation follows from (4.24), the fourth equation follows
from the definition of the function $Z_{M}^{(\alpha )}\ $and the fifth
equation follows the fourth equation after obvious manipulations. \ $%
\blacksquare $

\bigskip

Remark. When $V=\infty $, for $\alpha \geq 0\ $we let $\eta _{M}(\alpha
)=\sup \{\theta :\phi (\theta )-\theta M=\alpha \}$since $Z_{M}^{(\alpha
)\alpha }(x)=O(e^{\eta _{M}(\alpha )x})$ as $x\rightarrow \infty ,\ $then we
have

\begin{eqnarray*}
E_{x}[\overset{\ast }{T^{-}}_{\tau }] &=&\frac{(x-\tau )}{\eta _{M}(\alpha
)^{^{\prime }}(0)} \\
&=&\frac{(x-\tau )}{M-E(I_{1})}\ ,\ \text{if }M>E(I_{1}) \\
&=&\infty \ \ \ \ \ \ \ \ \ \ \ \ \ \ \ \ ,\ \text{if\ }M\leq E(I_{1})\text{.%
}
\end{eqnarray*}%
This is consistent with the well known fact about the busy period of the
M/G/1 queuing system. \ \ \ \ \ \ \ \ \ \ \ \ \ \ \ \ \ \ \ \ \ \ \ \ \ \ \
\ \ \ \ \ \ \ \ \ \ \ \ \ \ \ \ \ \ \ \ \ \ \ \ \ \ \ \ \ \ \ \ \ \ \ \ \ \
\ \ \ \ \ \ \ \ \ \ \ \ \ \ \ 

\bigskip \bigskip To compute $E_{x}[\exp (-\alpha \overset{\ast }{T}_{0})]$,
we first observe that, for $\lambda \leq x\leq V$, $\overset{\ast }{T}_{0}=%
\overset{\ast }{T^{-}}_{\tau }\ $almost everywhere. Hence $E_{x}[\exp
(-\alpha \overset{\ast }{T}_{0})]$ is given in (4.25). We now turn our
attention to the case where $x<\lambda $.\ We first consider the case where
the input process is a spectrally positive L%
\'{}%
evy process.

\bigskip Theorem. Assume that the input process is a spectrally positive L%
\'{}%
evy process. For $z>\lambda $, we define

\bigskip 
\begin{equation*}
h_{\alpha }(x,dz)=\int_{0}^{\lambda }dyu^{\alpha }(x,y)\upsilon (dz-y)\text{,%
}
\end{equation*}
$u^{\alpha }(x,y)$ is defined in (4.6).

Then, for $\alpha \geq 0,x<\lambda $

\begin{equation*}
E_{x}[e^{-\alpha \overset{\ast }{T}_{0}}]=\frac{1}{Z^{\alpha }(\lambda -\tau
)}[\int_{\lambda }^{V}Z^{\alpha }(\lambda -z)h_{\alpha
}(x,dz)+\int_{V}^{\infty }h_{\alpha }(x,dz)]\ \ \ \ (4.27)\ 
\end{equation*}

\bigskip Proof:\ We write\ \ \ \ \ \ \ \ \ \ \ \ \ \ \ \ \ \ \ \ \ \ \ \ \ \
\ \ \ \ \ \ \ \ \ \ \ \ \ \ \ \ \ \ \ \ \ \ \ \ \ \ \ \ \ \ \ \ \ \ \ \ \ \
\ \ \ \ \ \ \ \ \ \ \ \ \ \ \ \ \ \ \ \ \ \ \ \ \ \ \ \ \ \ \ \ \ \ \ \ \ \
\ \ \ \ \ \ \ \ \ \ \ \ \ \ \ \ \ \ \ \ \ \ \ \ \ \ \ \ \ \ \ \ \ \ \ \ \ \
\ \ \ \ \ \ \ \ \ \ \ \ \ \ \ \ \ \ \ \ \ \ \ \ \ \ \ \ \ \ \ \ \ \ \ \ \ \
\ \ \ \ \ \ \ \ \ \ \ \ \ \ \ \ 

\begin{eqnarray*}
E_{x}[e^{-\alpha \overset{\ast }{T}_{0}}] &=&E_{x}[e^{-\alpha T_{\lambda
}^{+}-\alpha (\overset{\ast }{T}_{0}-T_{\lambda }^{+})}] \\
&=&E_{x}[E_{x}[e^{-\alpha T_{\lambda }^{+}-\alpha (\overset{\ast }{T}%
_{0}-T_{\lambda }^{+})}\mid \sigma (T_{\lambda }^{+},I_{T_{\lambda }^{+}}]]
\\
&=&E_{x}[e^{-\alpha T_{\lambda }^{+}\ \ \ }E_{(I_{T_{\lambda }^{+}}\wedge
V)}[e^{-\alpha T_{\tau }^{-}}]] \\
&=&\frac{1}{Z^{\alpha }(V-\tau )}E_{x}[e^{-\alpha T_{\lambda }^{+}\ \ \
}Z_{M}^{\alpha }(V-(I_{T_{\lambda }^{+}}\wedge V))] \\
&=&\frac{1}{Z_{M}^{\alpha }(V-\tau )}[\int_{0}^{V}Z_{M}^{\alpha
}(V-z)h_{\alpha }(x,dz)+\int_{V}^{\infty }h_{\alpha }(x,dz)]\text{,}
\end{eqnarray*}%
where the third equation follows since, given $T_{\lambda }^{+}$ and $%
I_{T_{\lambda }^{+}}$, $\overset{\ast }{T}_{0}-T_{\lambda }^{+}\ $is equal
to $T_{\tau }^{-}\ $almost everywhere. The fourth equations from (4.11), the
last equation follows from the fact that $Z^{\alpha }(0)=1$, and the
definition of $h_{\alpha }(x,dz)$.$\ $ \ \ $\blacksquare $

\bigskip

The following theorem gives a result analogous to (4.27) when the input
process is a spectrally positive L%
\'{}%
evy process \ reflected at its infimum.

\bigskip Theorem. Assume that the input process is a spectrally positive L%
\'{}%
evy process\ reflected at its infimum. For $z\geq \lambda $, let $l_{\alpha
}(dz)$, $L_{\alpha }(z)$, and $V_{\alpha }(\lambda )$\ be as defined in
(4.16), (4.17) and (4.18), respectively. Define

\begin{equation*}
g_{\alpha }(x,dz)=\left\{ 
\begin{array}{c}
=\frac{l_{\alpha }(dz)}{W_{+}^{(\alpha )^{\prime }}(\lambda )}\ \text{,\ }%
z>\lambda \\ 
=\frac{V_{\alpha }(\lambda )-L_{\alpha }(\lambda )}{W_{+}^{(\alpha )^{\prime
}}(\lambda )}\delta _{\lambda }(dz)\text{.}%
\end{array}%
\right.
\end{equation*}

\ \ Then, for $\alpha \geq 0,x<\lambda $

\begin{equation*}
E_{x}[e^{-\alpha \overset{\ast }{T}_{0}}]=\frac{1}{Z_{M}^{\alpha }(\lambda
-\tau )}[\int_{0}^{V}Z_{M}^{\alpha }(\lambda -z)g_{\alpha
}(x,dz)+\int_{V}^{\infty }g_{\alpha }(x,dz)]\text{. \ \ \ \ }(4.28)\ \ \ 
\end{equation*}

Proof. The proof follows in a manner similar to the proof of (4.27), using
(4.19), (4.20) and (4.25). \ \ \ \ \ $\blacksquare $\ \ \ \ \ \ \ \ \ \ \ \
\ \ \ \ \ \ \ \ \ \ \ \ \ \ \ \ \ \ \ \ \ \ \ \ \ \ \ \ \ \ \ \ \ \ 

\bigskip \noindent {\LARGE 5.Special Cases}

\medskip {\LARGE \ }

In this section we consider the cases where the input process is a
spectrally positive L\'{e}vy \ of bounded variation, Brownian motion
reflected at its infimum and Wiener process. For the first case, we extend
the results of Bae\textit{\ et al} \ [5] who assumed that the input process
is a compound Poisson process with a negative drift. We also simplify some
of their results. For the second case, we obtain results similar to those of
Attia [3] and Lam and Lou [7]. In the third case we obtain results similar
to those of Zuckerman [10].

{\LARGE (i) }Assume that the input is a spectrally positive L\'{e}vy process
of bounded variation with L\'{e}vy exponent described in (3.3), reflected at
its infimum. Let \noindent\ $\mu =\int_{0}^{\infty }x\upsilon (dx)$ and
assume that\ $\mu <\infty $. For every $x\in R_{+}$, we define the
probability density function $f(x)=\frac{\upsilon ([x,\infty ))}{\mu }$, and 
$\int_{0}^{\infty }(1-e^{-\theta x})\upsilon (dx)=$ $\theta \mu
\int_{0}^{\infty }\ e^{-\theta x}f(x)dx$.\ Define $\rho =\frac{\mu }{%
\varsigma }$ and assume that $\rho <1$, then, $\frac{1}{\phi (\theta )}=%
\frac{1}{\varsigma \theta \lbrack 1-\rho \int_{0}^{\infty }\ e^{-\theta
x}f(x)dx]}=\frac{1}{\varsigma \theta }\int_{0}^{\infty }e^{-\theta
x}dx\sum\limits_{n=0}^{\infty }\rho ^{n}f^{(n)}(x)=\frac{1}{\varsigma }%
\int_{0}^{\infty }e^{-\theta x}dx\sum\limits_{n=0}^{\infty }\rho
^{n}F^{(n)}(x)$, where $F(x)\ $is the distribution function corresponding to 
$f\ $. Therefore, the $0-$scale function is given is given as follows

\begin{equation*}
W(x)=\frac{1}{\varsigma }\sum\limits_{n=0}^{\infty }\rho ^{n}F^{(n)}(x)\text{%
. \ \ }
\end{equation*}%
For $\alpha >0,W^{(\alpha )}$ is computed using $(3.5)\ $and the above
equation.

\bigskip

Remark. (a) Since, for all $n\geq 0$, $F^{(n)}(x)\leq \lbrack F(x)]^{n}$,
then for all $x$, $W(x)\leq \frac{1}{\varsigma -\mu F(x)}$,\ if $\rho <1$.

(b) Bae\textit{\ et al} \ [5], treat the special case where the input
process is a compound Poisson process with a negative drift. In this case, $%
\upsilon (dx)=\lambda G(dx)$,$\,\ $where $\lambda >0$ and $G\ $\ is a
distribution function of a positive random variable $[0,\infty )$,
describing the size of each jump of the compound Poisson process. $\ $In
this case, $f(x)=\frac{\overset{-}{G}(x)\ }{m}$ and $\rho =\frac{\lambda m}{%
\varsigma },\ $where $\overset{-}{G}=1-G$\ \ and $m$ $=\int_{0}^{\infty }%
\overset{-}{G}(x)dx$, which is assumed to be finite. We note that their
entities $w(x)$ and $E[L^{a}(\lambda ,\tau )]$ given in p. 521 are nothing
but our $C(x,\lambda ,0)$ and $E_{x}[\tau _{\lambda }\ ]$, respectively. $\ $%
Using \ (4.12) and (4.15) we get the same result given in page 523 of this
reference.\ Furthermore, their functions $E[P^{M}\ (\lambda ,\tau )]\ $and $%
E[L^{M}\ (\lambda ,\tau )]\ $given in page 525 are our $C(x,M,\tau )\ $and $%
E_{x}[\overset{\ast }{T^{-}}_{\tau }]$, respectively. Using (4.24) and
(4.26) we get a simpler form for these entities. Furthermore, letting $%
\alpha =0$, in (4.19) we provide a simpler formula for the distribution of $%
L(\tau )\ $(the overshoot) given on page 525 of their paper.

\bigskip (c) Assume that the input process is a gamma process with negative
drift. The L\'{e}vy measure is given by $\upsilon (dx)=a\frac{e^{-bx}}{x}dx$%
{\LARGE , }$a,b>0.\ $In this case, $E(I_{1})=\varsigma -\frac{a}{b}$, which
is assumed to be nonnegative and $\rho =\frac{a}{\varsigma b}$\ $<1$.$\ $It
follows that $f(x)=b\int_{x}^{\infty }${\LARGE \ }$\frac{e^{-by}{\LARGE \ }}{%
y}$, the \ right hand side is denoted by $E_{1}(x)\ $in p. 227 of [2], \
Direct integrations yields, $F(x)=(1-e^{-bx})+xf(x)$.$\ \ \ $\ \ \ \ 

\ (d) Assume that the input process is an inverse Gaussian process with a
negative drift, and with L\'{e}vy measure is given by $\upsilon (dx)=\frac{1%
}{\sigma \sqrt{2\pi x^{3}}}e^{-xc^{2}/2\sigma ^{2}}${\LARGE , }$\sigma ,c>0$%
. It follows that $E(I_{1})=\varsigma -\frac{1}{c}$, which is assumed to be
grater than zero. In this case $\rho =$ $\frac{1}{c\varsigma }<1$. In this
case, $f(x)=c\int_{x}^{\infty }{\LARGE \ }\upsilon (dy)$, and $F(x)=\func{erf%
}(c\sqrt{y/2\sigma ^{2}})+xf(x)$.\ \ \ \ \ \ \ \ \ \ \ \ \ \ \ \ \ \ \ \ \ \
\ \ \ \ \ \ \ 

.$\ \ \ \ \ \ \ \ \ \ \ \ \ \ \ \ \ \ \ \ \ \ \ \ \ \ \ \ \ \ \ \ \ \ \ \ \
\ \ \ \ \ \ \ \ \ \ \ \ \ \ \ \ \ \ \ \ \ \ \ \ \ \ \ \ \ \ \ \ \ \ \ \ \ \
\ \ \ \ $

\noindent {\LARGE (ii) \ }Assume that the input process is a Brownian motion
with drift term $\mu \in R$, variance term $\sigma ^{2}$, reflected at its
infimum. We will show that the results of [3] and [7] follow from our
results.\ In this case, the L\'{e}vy measure $\nu =0$, and from (3.1) we
have, that for $\theta \geq 0$, $\varphi (\theta )=-\mu \theta +\frac{\theta
^{2}\sigma ^{2}}{2}$. It follows that, for $\alpha \geq 0$, $\eta (\alpha )=%
\frac{\sqrt{2\alpha \sigma ^{2}+\mu ^{2}}+\mu }{\sigma ^{2}}$. Let $\delta =%
\sqrt{2\alpha \sigma ^{2}+\mu ^{2}}$, we have, $W^{\alpha }(x)=\frac{2}{%
\delta }e^{\mu x/\sigma ^{2}}\sinh (\frac{x\delta }{\sigma ^{2}})\ $and $\
Z^{\alpha }(x)=e^{\mu x/\sigma ^{2}}\left( \cosh (\frac{x\delta }{\sigma ^{2}%
})-\frac{\mu }{\delta }\sinh (\frac{x\delta }{\sigma ^{2}})\right) .\ \ $We
note that $W^{\alpha }(x)\ $is differentiable, and $W^{\alpha ^{\prime }}(x)=%
\frac{\mu }{\sigma ^{2}}W^{\alpha }(x)+\frac{2}{\sigma ^{2}}e^{\mu x/\sigma
^{2}}\cosh (\frac{x\delta }{\sigma ^{2}})$, it follows that $\frac{%
W^{(\alpha )}(\lambda )}{W^{(\alpha )^{\prime }}(\lambda )}=\left( \frac{%
\sigma ^{2}}{\mu +\delta \coth ((\frac{\lambda \delta }{\sigma ^{2}})}%
\right) $. Substituting the values of $Z^{(\alpha )}(\lambda -x),W^{(\alpha
)}(\lambda -x)$ and $\frac{W^{(\alpha )}(\lambda )}{W^{(\alpha )^{\prime
}}(\lambda )}\ $in (4.14), we have, for $\alpha \geq 0,x\leq \lambda $

\begin{equation*}
E_{x}[e^{-\alpha {\huge \tau }_{_{\lambda }}}]=e^{\mu (\lambda -x)}\left[
\cosh \left( \frac{(\lambda -x)\delta }{\sigma ^{2}}\right) -\frac{1}{\delta 
}\sinh \left( \frac{(\lambda -x)\delta }{\sigma ^{2}}\right) \left( \mu +%
\frac{2\alpha \sigma ^{2}}{\mu +\delta \coth (\frac{\lambda \delta }{\sigma
^{2}})}\right) \right] \text{.}
\end{equation*}%
Case 1. $\mu \#0$: It follows that, for $x\geq 0,$ $W(x)=\frac{e^{2\mu
x/\sigma ^{2}}-1}{\mu },W^{^{\prime }}(x)=\frac{2e^{2\mu x/\sigma ^{2}}}{%
\sigma ^{2}}$ and $\overset{-}{W}(x)=\frac{\sigma ^{2}}{2\mu ^{2}}(e^{2\mu
x/\sigma ^{2}}-1)-\frac{x}{\mu }.\ $Substituting the values of $W(\lambda
-x) $, $\overset{-}{W}(\lambda -x),W(\lambda )\ $and $W^{^{\prime }}(\lambda
)\ $in (4.15) we have, for $x\leq \lambda $,

\begin{equation*}
E_{x}[{\huge \tau }_{_{\lambda }}]=\frac{\lambda -x}{\mu }+\frac{\sigma ^{2}%
}{2\mu ^{2}}\left[ e^{-2\mu \lambda /\sigma ^{2}}-e^{-2\mu x/\sigma ^{2}}%
\right] .
\end{equation*}

We note that, $W^{(\alpha )}(0)=0$, $\overset{(2)}{U^{\alpha }}$\ in (4.12)$%
\ $is absolutely continuous \ with respect to the Lebesgue measure on $%
[0,\lambda )\ $and for $y\in $ $[0,\lambda )$, $W^{(\alpha )}(dy)=W^{\alpha
^{\prime }}(y)dy$. Substituting the values of $W^{(\alpha )}(\lambda -x)$,\ $%
W^{\alpha ^{\prime }}(\lambda )$, $W^{(\alpha )}(y-x)$, and $W^{\alpha
^{\prime }}(y)\ $in (4.12) we get a version of the the density of $\overset{%
(2)}{U^{\alpha }}$,$\ $denoted by\ $\overset{(2)}{u^{\alpha }}$. Thus , $%
C_{\alpha }(x,\lambda ,0)$ is computed using (4.3).

Let$\ \overset{\ast }{\mu }^{{}}=\mu -M$,$\ \overset{\ast }{\delta }=\sqrt{%
2\alpha \sigma ^{2}+\overset{\ast }{\mu }^{2}}$, \ we note that $%
W_{M}^{\alpha }(x)=\frac{2}{\delta }e^{\overset{\ast }{\mu }x/\sigma
^{2}}\sinh (\frac{x\overset{\ast }{\delta }}{\sigma ^{2}})$.\ \ We note that
the input process is continuous, $Y_{_{T_{\lambda }^{+}}}=\lambda <V$, $\ $%
almost every where. Therefore, the term $C_{\alpha }((Z_{\overset{\wedge }{T}%
_{0}}\wedge V),\tau ,M)$ in (2.6) reduces to $C_{\alpha }(\lambda ,\tau ,M)$
which is computed using (4.23). Furthermore,

\begin{equation*}
E_{x}[e^{-\alpha \overset{\ast }{T}_{0}}]=E_{x}[e^{-\alpha {\huge \tau }%
_{_{\lambda }}}]E_{\lambda }[e^{-\alpha T_{\tau }^{-}}]\text{, }
\end{equation*}%
where $E_{x}[e^{-\alpha {\huge \tau }_{_{\lambda }}}]\ $is given above and $%
E_{\lambda }[e^{-\alpha T_{\tau }^{-}}]$ is given in (4.25).

Let $\lambda ^{\ast }=V-\lambda \ $and $\tau ^{\ast }=V-\tau $, then

\begin{equation*}
E_{\tau }[\overset{\ast }{T}_{0}]=E_{\tau }[{\huge \tau }_{_{\lambda
}}]+E_{\lambda }[T_{\tau }^{-}],
\end{equation*}%
Note that

\begin{equation*}
E_{\lambda }[T_{\tau }^{-}]=\frac{\lambda ^{\ast }-\tau ^{\ast }}{\mu ^{\ast
}}+\frac{\sigma ^{2}}{2\overset{\ast }{\mu }^{2}}\left[ e^{2\overset{\ast }{%
\mu }\tau ^{\ast }/\sigma ^{2}}-e^{2\mu \lambda ^{\ast }/\sigma ^{2}}\right] 
\text{,}
\end{equation*}%
where the last equation follows from (4.26) after some tedious calculations
which we omit.

\bigskip

Case 2. $\mu =0$: In this case, $\delta =\sqrt{2\alpha \sigma ^{2}}$,$\ \ $%
letting\ $\mu \rightarrow 0$, in the corresponding equations in case 1
above, we have

\begin{equation*}
E_{x}[e^{-\alpha T_{\lambda }^{+}}]=\left[ \cosh ((\lambda -x)\delta /\sigma
^{2})-\sigma ^{2}\frac{\sinh \left( (\lambda -x)\delta /\sigma ^{2}\right) }{%
\coth (\frac{\lambda \delta }{\sigma ^{2}})}\right] \text{,}
\end{equation*}

\begin{equation*}
E_{x}[T_{\lambda }^{+}]=\frac{\lambda ^{2}-x^{2}}{\sigma ^{2}}\text{,}
\end{equation*}%
and $E_{\lambda }[T_{\tau }^{-}]\ $is obtained by replacing $\mu ^{\ast }$
by $-M\ \ $in the last equation of case 1.

\bigskip

\noindent {\LARGE (iii)} Assume that the input process is a Brownian motion
with drift term $\mu \ >0$\ and variance parameter $\sigma ^{2}$.
Substituting the values of $W^{(\alpha )}(x),Z^{(\alpha )}(x)$, given in 
{\large (ii)} in (4.7) we have, for $x\leq \lambda $, $E_{x}[e^{-\alpha
T_{\lambda }^{+}}]=\exp \left( (\delta -\mu )(x-\lambda )\right) $.
Substituting$\ \eta (0)=\frac{2\mu }{\sigma ^{2}}$, $\frac{1}{\mu }(e^{2\mu
x/\sigma ^{2}}-1)$ and $\frac{\sigma ^{2}}{2\mu ^{2}}(e^{2\mu x/\sigma
^{2}}-1)-\frac{x}{\mu }$ for $W(x)\ $and $\overset{-}{W}(x)$, respectively,
in (4.8) we have, for $x\leq \lambda $, $E_{x}[T_{\lambda }^{+}]=\frac{%
\lambda -x}{\mu }$. These results are consistent with the results of
Zuckerman [10], p.423. The computations of the other entities in the cost
functionals (2.7)\ and (2.8) can be obtained in a manner similar to those
discusses in {\large (ii)} with obvious modifications.

\bigskip

{\LARGE References}

[1] Abdel-Hameed, M. (2000). Optimal control of a dam using $P_{\lambda
,\tau }^{M}\ $policies and penalty cost when the input process is a compound
Poisson process with positive drift. \textit{J.Appl.Prob. }\textbf{37},
508-416.

[2] Abramowitz, M. and Stegun, I. A. (1964). \textit{Handbook of
Mathematical Functions}. Dover, New York.

[3] Attia, F. (1987). The control of a finite dam with penalty cost
function; Wiener process input. \textit{Stochastic Processes and Their
Applications} \textbf{25},

289-299.

[4] Bae, J, Kim, S. and Lee, E.Y. (2002). A\ $P_{\lambda }^{M}$\ policy for
an M/G/1 queueing system. \textit{Appl. Math. Modelling} \ \textbf{26, }%
929-939.

[5] Bae, J, Kim, S. and Lee, E.Y. (2003). Average cost under $P_{\lambda
,\tau }^{M}\ $- policy in a finite dam with compound Poisson input. \textit{%
J.Appl.Prob. }\textbf{40}, 519-526.

[6] Kyprianou, A. E. (2006). \textit{Introductory Lecture Notes on
Fluctuations of L%
\'{}%
evy Processes with Applications}. Springer Verlag.

[7] Lam, Y. and Lou, J.H. (1987). Optimal control of a finite dam: Wiener
process input. \textit{J.Appl.Prob.} \textbf{35}, 482-488.

[8] Lea, E.Y. and Ahn, S.K. (1998). $P_{\lambda }^{M}$\ policy for a dam
with input formed by a compound Poisson process. \textit{J.Appl.Prob.} 
\textbf{24}, 186-199.

[9] Zhou, X. W. (2004). Some fluctuation identities for L%
\'{}%
evy processes with jumps of the same sign. \textit{J.Appl.Prob.} \textbf{41}%
, 1191-1198.

[10] Zuckerman, D. (1977). Two-stage output procedure for a finite dam. 
\textit{J.Appl.Prob.} \textbf{14}, 421-425.

\bigskip\ \ \ \ 

{\LARGE \ }

\end{document}